\documentclass[a4paper, 11pt]{amsart}

\usepackage{amsmath, amssymb, amsthm, mathtools}
\mathtoolsset{showonlyrefs=true}
\usepackage[colorlinks=true, linkcolor=blue]{hyperref}

\usepackage[margin=1in]{geometry}
\usepackage{quiver}
\usepackage{booktabs}
\usepackage{array}
\usepackage{float}

\newtheorem{theorem}{Theorem}[section]
\newtheorem{lemma}[theorem]{Lemma}
\newtheorem{proposition}[theorem]{Proposition}
\newtheorem{claim}[theorem]{Claim}

\newtheorem*{main}{Main Theorem}

\theoremstyle{definition}
\newtheorem{definition}[theorem]{Definition}

\newtheorem{remark}[theorem]{Remark}

\newtheorem{problem}[theorem]{Problem}

\newenvironment{pfclaim}
    {\begin{proof}[Proof of Claim]}
    {\end{proof}}

\newcommand{\PP}{\mathbb{P}}

\newcommand{\Spec}{\operatorname{Spec}}

\newcommand{\rad}{\operatorname{rad}}
\newcommand{\depth}{\operatorname{depth}}

\let\subset\subseteq

\renewcommand{\hat}[1]{\widehat{#1}}

\renewcommand{\tilde}[1]{\widetilde{#1}}

\title[A CM Analogue of the Brodmann--Rotthaus Theorem]{A Cohen--Macaulay Analogue of the Brodmann--Rotthaus Theorem}
\author{Rirai Ikeda}
\address{Department of Mathematics, School of Science, Science Tokyo, 2-12-1 Ookayama, Meguro-ku, Tokyo 152-8551, Japan}
\email{ikeda.r.cfbf@m.isct.ac.jp}

\begin{document}

\begin{abstract}
    We prove a Cohen--Macaulay analogue of the theorem of Brodmann and Rotthaus on the lifting of quasi-excellence. More precisely, let $A$ be a Noetherian ring and let $I$ be an ideal contained in the Jacobson radical of $A$. We show that if $A/I$ is $\mathsf{CM}$-Q0 and $A$ is a $\mathsf{CM}$\textnormal{-G-}ring, then $A$ is $\mathsf{CM}$-quasi-excellent. Our proof follows the strategy of Kurano and Shimomoto and makes essential use of \v{C}esnavi\v{c}ius' Macaulayfication theorem.
\end{abstract}

\maketitle

\tableofcontents

\section{Introduction}\label{sec:introduction}

Throughout the present paper, all rings are commutative with identity. For a ring $A$, we denote its Jacobson radical by $\rad(A)$.

The following problem, originating in a question of Grothendieck \cite[({\textbf{IV}},~7.4.8)]{Gro65}, is known as the ``lifting problem.''

\begin{problem}[Lifting problem]\label{prob:lifting}
    Let $A$ be a Noetherian ring, let $I \subset A$ be an ideal, and let $\PP$ be a property of rings. Assume that the following conditions hold:
    \begin{enumerate}
        \item $A$ is $I$-adically complete.
        \item $A/I$ satisfies $\PP$.
    \end{enumerate}
    Does $A$ also satisfy $\PP$?
\end{problem}

The lifting problem is now settled for many properties $\PP$; we refer the reader to \cite{KS21,Lyu25,Lyu26} for further details. 

In their study of this problem, Brodmann and Rotthaus proved the following result. 

\begin{theorem}[{\cite[Theorem~2]{BR80}}]
    Let $A$ be a Noetherian universally catenary ring containing a field of characteristic $0$. Suppose that $A$ is a $G$-ring and that $A/I$ is \textup{J-2}, where $I$ is an ideal contained in the Jacobson radical of $A$. Then $A$ is also \textup{J-2}.
\end{theorem}

The proof of this theorem relies on Hironaka's resolution of singularities \cite{Hir64}, which accounts for the assumption that $A$ contains a field of characteristic $0$. Indeed, Brodmann and Rotthaus remarked in \cite[p.~81]{BR80} that this assumption could be removed if resolution of singularities were available in arbitrary characteristic.

Using his weak local uniformization theorem \cite[Expos\'e~VII, Th\'eor\`eme~1.1]{ILO14}, Gabber later removed this assumption. His argument was outlined in a letter to Laszlo \cite{Gabber2007}, and Kurano and Shimomoto gave a complete proof in \cite{KS21}. Thus, we have the following characteristic-free version of the theorem of Brodmann and Rotthaus.

\begin{theorem}[{\cite{Gabber2007}, \cite[Theorem~4.1]{KS21}}]
    Let $A$ be a Noetherian ring with an ideal $I \subseteq \operatorname{rad}(A)$. Assume that $A/I$ is quasi-excellent and that $A$ is a $G$-ring. Then $A$ is a \textup{J-2} ring. In other words, $A$ is quasi-excellent.
\end{theorem}

In fact, the assumption that $A/I$ is quasi-excellent in the above theorem can be weakened to the assumption that $A/I$ satisfies J-2.

The proof of the theorem of Brodmann and Rotthaus discussed above relies on resolution of singularities. Its Cohen--Macaulay version is Macaulayfication, in which regularity is replaced by the Cohen--Macaulay property. This suggests that Macaulayfication may provide a similar approach to the Cohen--Macaulay setting. The main theorem of the present paper confirms this expectation. 

\begin{main}[Theorem~\ref{thm:main1}]
    Let $A$ be a Noetherian ring with an ideal $I \subset \rad(A)$. Assume that $A/I$ is $\mathsf{CM}$\textup{-Q0} and that $A$ is a $\mathsf{CM}$\textnormal{-G-}ring. Then $A$ is a $\mathsf{CM}$\textnormal{-Q0} ring. In other words, $A$ is $\mathsf{CM}$-quasi-excellent.
\end{main}

The $\mathsf{CM}$-Q0 condition is a Cohen--Macaulay analogue of J-2; see Definition~\ref{def:CM-quasi-excellent} and Remark~\ref{rem:CM-quasi-excellent} for details.

Recall that a Macaulayfication of a Noetherian scheme $X$ is a projective birational morphism $\widetilde{X} \to X$ such that $\widetilde{X}$ is Cohen--Macaulay. The Macaulayfication problem was initiated by Faltings and subsequently studied by several authors. In particular, Kawasaki established a series of Macaulayfication results in \cite{Kaw98,Kaw00,Kaw02,Kaw08}; see \cite[\S~1.12]{Ces21} for a detailed account of the previous work. For our purposes, \v{C}esnavi\v{c}ius' theorem is particularly well suited, since it establishes the existence of a Macaulayfication under the sole assumption of $\mathsf{CM}$-quasi-excellence \cite[Theorem~1.6]{Ces21}. We use this theorem in Lemma~\ref{lem:lem1-3}, which is a key step in the proof of the main theorem.

The relationship between the results discussed above and the methods used in their proofs is summarized in the following diagram.

\begin{center}
  \begin{tikzcd}[
    column sep=-6pt,
    row sep=normal
  ]
    {\text{Regular case over }\mathbb{Q}} & {\text{Regular case}} & {\text{CM case}} \\[-1.5em]
    \begin{array}{c}
      \text{Resolution of singularities} \\
      \text{\cite{Hir64}}
    \end{array} &
    \begin{array}{c}
      \text{Resolution of singularities} \\
      \text{(open problem)}
    \end{array} &
    \begin{array}{c}
      \text{Macaulayfication} \\
      \text{\cite{Ces21}}
    \end{array} \\
    &
    \begin{array}{c}
      \text{Gabber's weak local uniformization} \\
      \text{\cite[Expos\'e~VII, Th\'eor\`eme~1.1]{ILO14}}
    \end{array} &
    \\
    \begin{array}{c}
      \text{Brodmann--Rotthaus theorem} \\
      \text{\cite[Theorem~2]{BR80}}
    \end{array} &
    \begin{array}{c}
      \text{Gabber--Kurano--Shimomoto theorem} \\
      \text{\cite{Gabber2007}, \cite[Theorem~4.1]{KS21}}
    \end{array} &
    \begin{array}{c}
      \text{Main theorem} \\
      \text{Theorem~\ref{thm:main1}}
    \end{array}
    \arrow[ Rightarrow, from=2-1, to=4-1 ]
    \arrow[ Rightarrow, from=2-2, to=3-2 ]
    \arrow[ Rightarrow, from=2-3, to=4-3 ]
    \arrow[ Rightarrow, from=3-2, to=4-2 ]
  \end{tikzcd}
\end{center}

Our proof follows the strategy of Kurano and Shimomoto. The use of Macaulayfication simplifies the part of their argument involving alteration coverings. On the other hand, unlike in a regular local ring, the maximal ideal of a Cohen--Macaulay local ring need not be generated by a regular sequence. For this reason, the corresponding part of their argument requires some modification; see Lemma~\ref{lem:lem1-1}.

The paper is organized as follows. In Section~\ref{sec:preliminaries}, we recall the notions of $\mathsf{CM}$-G-rings, $\mathsf{CM}$-Q0 rings, and $\mathsf{CM}$-quasi-excellent rings and establish several preliminary results. In Section~\ref{sec:CM-version}, we prove the main theorem.

\medskip

\noindent
\textbf{Acknowledgements.} The author would like to thank his supervisor, Professor Kazuma Shimomoto, for his helpful comments and suggestions. The author is grateful to Kazuki Hayashi for his valuable comments on the first draft of the paper.

\section{Preliminaries for Cohen--Macaulay analogues}\label{sec:preliminaries}

In this section, we show that several basic properties of quasi-excellent rings have analogues for $\mathsf{CM}$-quasi-excellent rings. From now on, for a scheme $X$ and a ring $A$, we denote their Cohen--Macaulay loci by $\mathsf{CM}(X)$ and $\mathsf{CM}(A)$, respectively.

Recall that the definition of quasi-excellence consists of two conditions:
\begin{enumerate}
    \item[(G-ring)] A Noetherian ring $A$ is a \emph{G-ring} if for every prime ideal $\mathfrak{p}$ of $A$, all formal fibers $\hat{A_\mathfrak{p}} \otimes_{A_\mathfrak{p}} \kappa(\mathfrak{q})$ of $A_\mathfrak{p}$ are geometrically regular over $\kappa(\mathfrak{q})$, where $\mathfrak q$ is a prime ideal of $A_\mathfrak p$ and $\kappa(\mathfrak q)$ is its residue field.
    \item[(J-2)] A Noetherian ring $A$ is \emph{J-2} if for any finitely generated $A$-algebra $B$, the regular locus $\mathsf{Reg}(B)$ is open in $\Spec B$. 
\end{enumerate}
Moreover, a Noetherian ring $A$ is \emph{excellent} if it is quasi-excellent and universally catenary. Replacing regularity in these conditions by Cohen–Macaulayness leads to the following definitions.

\begin{definition}\label{def:CM-quasi-excellent}
    Let $A$ be a Noetherian ring. 
    \begin{enumerate}
        \item We say that $A$ is a \emph{$\mathsf{CM}$\textnormal{-G-}ring} if for every prime ideal $\mathfrak{p}$ of $A$, all formal fibers of $A_\mathfrak{p}$ are Cohen--Macaulay.
        \item We say that $A$ is a \emph{$\mathsf{CM}$\textnormal{-Q0} ring} if for every prime ideal $\mathfrak{p}$ of $A$, the set $\mathsf{CM}(A/\mathfrak{p})$ contains a non-empty open subset of $\Spec(A/\mathfrak{p})$.
    \end{enumerate}
    We say that $A$ is \emph{$\mathsf{CM}$-quasi-excellent} if $A$ is a $\mathsf{CM}$\textnormal{-G-}ring and a $\mathsf{CM}$\textnormal{-Q0} ring. Moreover, we say that $A$ is \emph{$\mathsf{CM}$-excellent} if $A$ is $\mathsf{CM}$-quasi-excellent and universally catenary.
\end{definition}

For affine schemes, this definition agrees with the definition of $\mathsf{CM}$-quasi-excellence and $\mathsf{CM}$-excellence given by \v{C}esnavi\v{c}ius in \cite[Definition~1.2]{Ces21}. Clearly, any quasi-excellent ring is $\mathsf{CM}$-quasi-excellent. 

\begin{remark}\label{rem:CM-quasi-excellent}
    The notion of $\mathsf{CM}$-quasi-excellence is a Cohen--Macaulay version of the notion of quasi-excellence. Indeed, over a field, Cohen--Macaulayness is equivalent to geometric Cohen--Macaulayness. More precisely, a Noetherian ring $A$ over a field $k$ is Cohen--Macaulay if and only if $A \otimes_k L$ is Cohen--Macaulay for any finite field extension $L/k$. This is a consequence of \cite[\href{https://stacks.math.columbia.edu/tag/045N}{Tag 045N}]{stacks-project} or \cite[Theorem~23.6 and its Remark]{Mat86}. Moreover, $\mathsf{CM}$-Q0 and the following $\mathsf{CM}$-2 are equivalent; see \cite[Proposition 3.10]{Ike26}, applied with $\mathbb{P} = \mathsf{CM}$. 
    \begin{enumerate}
        \item[($\mathsf{CM}$-2)] A Noetherian ring $A$ is a $\mathsf{CM}$\textnormal{-2} ring if for every finitely generated $A$-algebra $B$, the Cohen--Macaulay locus $\mathsf{CM}(B)$ is open in $\Spec B$. 
    \end{enumerate}
\end{remark}

A $\mathsf{CM}$-G-ring is precisely a $P$-ring in the sense of \cite[\href{https://stacks.math.columbia.edu/tag/0BIR}{Tag 0BIR}]{stacks-project}, where $P(k \to R)$ means that $R$ is Cohen--Macaulay. Hence, we have the following basic properties. 

\begin{proposition}\label{prop:CM-G}
    Let $A$ be a Noetherian ring. Then the following statements hold:
    \begin{enumerate}
        \item If for every maximal ideal $\mathfrak{m}$ of $A$, $A_\mathfrak{m}$ is a $\mathsf{CM}$\textnormal{-G-}ring, then $A$ is a $\mathsf{CM}$\textnormal{-G-}ring; see \cite[\href{https://stacks.math.columbia.edu/tag/0BIU}{Tag 0BIU}]{stacks-project}. 
        \item If $A$ is a $\mathsf{CM}$\textnormal{-G-}ring, then every $A$-algebra essentially of finite type is also a $\mathsf{CM}$\textnormal{-G-}ring; see \cite[\href{https://stacks.math.columbia.edu/tag/0BIV}{Tag 0BIV}]{stacks-project}. 
        \item Let $A \to B$ be a ring homomorphism between Noetherian rings, and let $A'$ be an $A$-algebra of finite type. If all fibers of $A \to B$ are Cohen--Macaulay, then all fibers of the base change $A' \to A' \otimes_A B$ are Cohen--Macaulay; see \cite[\href{https://stacks.math.columbia.edu/tag/0BIR}{Tag 0BIR} (A)]{stacks-project}.
    \end{enumerate}
\end{proposition}

Every semi-local G-ring is J-2. The analogous property holds for $\mathsf{CM}$-G-rings. 

\begin{proposition}\label{prop:semilocalCMG}
    A semi-local $\mathsf{CM}$\textnormal{-G-}ring is also a $\mathsf{CM}$\textnormal{-Q0} ring. In particular, a semi-local $\mathsf{CM}$\textnormal{-G-}ring is $\mathsf{CM}$-quasi-excellent.
\end{proposition}
\begin{proof}
    First, we show the following claim. 
    \begin{claim}\label{claim:claimCMG}
        Let $A, B$ be Noetherian rings, and let $f \colon A \to B$ be a faithfully flat morphism whose fibers are Cohen--Macaulay. If $B$ is a $\mathsf{CM}$\textnormal{-Q0} ring, then so is $A$.
    \end{claim}
    \begin{pfclaim}
        Note that $\mathsf{CM}$-Q0 and $\mathsf{CM}$-2 are equivalent. Therefore, taking base change to a finite type $A$-algebra, it is enough to show that $\mathsf{CM}(A)$ is Zariski open in $\Spec A$ by Proposition \ref{prop:CM-G} (3).

        Since $f$ is faithfully flat, the induced map $f^* \colon \Spec B \to \Spec A$ is submersive. By \cite[Corollary (21.E) iii)]{Mat70}, we have $(f^*)^{-1}(\mathsf{CM}(A)) = \mathsf{CM}(B)$. Hence, $\mathsf{CM}(A)$ is open in $\Spec A$.
    \end{pfclaim}
    We now show Proposition \ref{prop:semilocalCMG}. Let $A$ be a semi-local $\mathsf{CM}$-G-ring. Let $A \to \hat{A}$ be the completion of $A$ with respect to the Jacobson radical. We have
    \begin{align*}
        \hat{A} \cong \prod_{i=1}^n \hat{A_{\mathfrak{m}_i}},
    \end{align*}
    where $\mathfrak{m}_1, \ldots, \mathfrak{m}_n$ are the maximal ideals of $A$. Recall that $A \to \hat{A}$ is faithfully flat. Since $A$ is a $\mathsf{CM}$-G-ring, all fibers are Cohen--Macaulay. Moreover, because complete local rings are $\mathsf{CM}$-excellent, so is $\hat{A}$. By Claim \ref{claim:claimCMG}, $A$ is a $\mathsf{CM}$-Q0 ring. This completes the proof of Proposition \ref{prop:semilocalCMG}. 
\end{proof}

The following proposition is an analogue of \cite[Theorem 3.1]{Gre76}. 

\begin{proposition}\label{prop:greco}
    Let $A \to A'$ be a module-finite ring homomorphism which induces a surjective map on the spectra. Then the following statements hold:
    \begin{enumerate}
        \item $A$ is $\mathsf{CM}$-quasi-excellent if and only if $A'$ is $\mathsf{CM}$-quasi-excellent. 
        \item $A$ is $\mathsf{CM}$-excellent if and only if it is universally catenary and $A'$ is $\mathsf{CM}$-quasi-excellent.
    \end{enumerate}
\end{proposition}
\begin{proof}
    (1) Note that we may replace ``J2'' by ``$\mathsf{CM}$-Q0'' in \cite[Proposition 2.4]{Gre76}. Moreover, in \cite[Proposition~1.3]{Gre76}, we can take ``$P$'' to be Cohen--Macaulay. Hence, we can prove Proposition~\ref{prop:greco} by the same argument as in the quasi-excellent case.

    (2) This is immediate from (1). 
\end{proof}

\section{A Cohen--Macaulay Analogue of the Brodmann--Rotthaus Theorem}\label{sec:CM-version}

In the present section, we prove the main theorem of the present paper, which is a Cohen--Macaulay analogue of \cite[Theorem~2]{BR80} and \cite[Theorem~4.1]{KS21}.

\begin{theorem}\label{thm:main1}
    Let $A$ be a Noetherian ring with an ideal $I \subset \rad(A)$. Assume that $A/I$ is $\mathsf{CM}$\textup{-Q0} and that $A$ is a $\mathsf{CM}$\textnormal{-G-}ring. Then $A$ is a $\mathsf{CM}$\textnormal{-Q0} ring. In other words, $A$ is $\mathsf{CM}$-quasi-excellent.
\end{theorem}

In order to show Theorem~\ref{thm:main1}, we prepare some lemmas. 

The following lemma already appeared in \cite[Lemma 2.1]{Tak07}, but we give a proof for the sake of completeness.

\begin{lemma}[{\cite[Lemma 2.1]{Tak07}}]\label{lem:lem1-1}
    Let $A$ be a Noetherian ring and $\mathfrak{p} \in \Spec A$ a prime ideal. Assume the following two conditions:
    \begin{enumerate}
        \item $A_\mathfrak{p}$ is Cohen--Macaulay. 
        \item $\mathsf{CM}(A/\mathfrak{p})$ contains a nonempty open subset of $\Spec A/\mathfrak{p}$. 
    \end{enumerate}
    Then, there exists an element $f \in A \setminus \mathfrak{p}$ such that $V(\mathfrak{p}) \cap D(f) \subset \mathsf{CM}(A)$. 
\end{lemma}
\begin{proof}
    Set $d = \dim A_\mathfrak{p}$. Since $A_\mathfrak{p}$ is Cohen--Macaulay, there exists a sequence $x_1, \ldots, x_d \in \mathfrak{p}$ which forms an $A_\mathfrak{p}$-regular sequence and $\sqrt{(x_1, \ldots, x_d)A_\mathfrak{p}} = \mathfrak{p}A_\mathfrak{p}$. Let $J = (x_1, \ldots, x_d)$. We can choose an element $g_1 \in A \setminus \mathfrak{p}$ such that the following holds:
    \begin{align*}
        \sqrt{JA[g_1^{-1}]} = \mathfrak{p}A[g_1^{-1}]. 
    \end{align*}
    Moreover, taking a suitable element $g_2 \in A \setminus \mathfrak{p}$, we have: 
    \begin{align*}
        \text{The sequence } x_1,\ldots, x_d \text{ forms an } A[g_2^{-1}] \text{-regular sequence in } A[g_2^{-1}]. 
    \end{align*}    
    On the other hand, we have the filtration:
    \begin{align*}
        A'/J' \supset (\mathfrak{p}' + J')/J' \supset (\mathfrak{p}'^2 + J')/J' \supset \cdots \supset (\mathfrak{p}'^N + J')/J' = 0, 
    \end{align*}
    where $A' = A[g_1^{-1}]$, $J' = JA[g_1^{-1}]$ and $\mathfrak{p}' = \mathfrak{p}A[g_1^{-1}]$. Let $N \geq 1$ be the minimum integer such that $(\mathfrak{p}'^N + J')/J' = 0$. For $i = 0, 1, \ldots, N$, we put $M_i = (\mathfrak{p}'^i + J')/J'$. Every $M_i/M_{i+1}$ for $i = 0, \ldots, N-1$ can be regarded as a finitely generated $A'/\mathfrak{p}'$-module naturally. By generic freeness, we can take an element $g_3 \in A \setminus \mathfrak{p}$ such that $M_i/M_{i+1} \otimes_{A'} A'[g_3^{-1}]$ is a free $A'[g_3^{-1}]/\mathfrak{p}'A'[g_3^{-1}]$-module for every $i = 0, \ldots, N-1$. Moreover, by our assumption (2), we can take an element $g_4 \in A \setminus \mathfrak{p}$ such that $A[g_4^{-1}]/\mathfrak{p}A[g_4^{-1}]$ is Cohen--Macaulay. Hence, setting $f = g_1g_2g_3g_4$, we have:
    \begin{enumerate}
        \item $x_1, \ldots, x_d$ form an $A[f^{-1}]$-regular sequence in $A[f^{-1}]$,
        \item $\sqrt{(x_1, \ldots, x_d)A[f^{-1}]} = \mathfrak{p}A[f^{-1}]$,
        \item $(M_i/M_{i+1}) \otimes_A A[f^{-1}]$ is a free $A[f^{-1}]/\mathfrak{p}A[f^{-1}]$-module $(\forall i = 0, \ldots, N-1)$, and 
        \item $A[f^{-1}]/\mathfrak{p}A[f^{-1}]$ is Cohen--Macaulay. 
    \end{enumerate}

    We now show that $V(\mathfrak{p}) \cap D(f) \subset \mathsf{CM}(A)$. Let $\mathfrak{q} \in V(\mathfrak{p}) \cap D(f)$. We put $t = \dim A_\mathfrak{q}/\mathfrak{p}A_\mathfrak{q}$. We have:
    \begin{align*}
        \depth_{A_\mathfrak{q}}((M_i/M_{i+1}) \otimes_A A_\mathfrak{q}) = \depth_{A_\mathfrak{q}/\mathfrak{p}A_\mathfrak{q}}((M_i/M_{i+1}) \otimes_A A_\mathfrak{q}) \geq t
    \end{align*}
    for every $i = 0, \ldots, N-1$, because $(M_i/M_{i+1}) \otimes_A A_\mathfrak{q}$ is a free $A_\mathfrak{q}/\mathfrak{p}A_\mathfrak{q}$-module and $A_\mathfrak{q}/\mathfrak{p}A_\mathfrak{q}$ is Cohen--Macaulay. (Note that $\depth 0 = \infty$.) For $i = 0, \ldots, N-1$, we have:
    \begin{align*}
        0 \to (M_{i+1})_\mathfrak{q} \to (M_i)_\mathfrak{q} \to (M_i/M_{i+1})_\mathfrak{q} \to 0. 
    \end{align*}
    We apply the depth lemma \cite[\href{https://stacks.math.columbia.edu/tag/00LX}{Tag 00LX}]{stacks-project} to the above exact sequence repeatedly from $i = N-1$ to $i = 0$. Then, we have:
    \begin{align*}
        t \leq \depth_{A_\mathfrak{q}}((M_0)_\mathfrak{q}) = \depth_{A_\mathfrak{q}}(A_\mathfrak{q}/JA_\mathfrak{q}) = \depth_{A_\mathfrak{q}/JA_\mathfrak{q}}(A_\mathfrak{q}/JA_\mathfrak{q}) \leq \dim A_\mathfrak{q}/JA_\mathfrak{q} = t.
    \end{align*}
    Hence, $A_\mathfrak{q}/JA_\mathfrak{q}$ is Cohen--Macaulay. Since $JA_\mathfrak{q}$ is generated by an $A_\mathfrak{q}$-regular sequence $x_1, \ldots, x_d$ in $A_\mathfrak{q}$, we have that $A_\mathfrak{q}$ is Cohen--Macaulay. Therefore, we have $V(\mathfrak{p}) \cap D(f) \subset \mathsf{CM}(A)$. This completes the proof of Lemma~\ref{lem:lem1-1}.
\end{proof}

Next, we show the following lemma, which is a Cohen--Macaulay analogue of \cite[Lemma~4.2]{KS21}. 

\begin{lemma}\label{lem:lem1-2}
    Let $A$ be a Noetherian ring, $I$ an ideal of $A$ and $\pi \colon X \to \Spec A$ a morphism of schemes which is of finite type. Assume that $A/I$ is $\mathsf{CM}$\textup{-Q0}. Then $\mathsf{CM}(X) \cap \pi^{-1}(V(I))$ is open in $\pi^{-1}(V(I))$.
\end{lemma}
\begin{proof}
    We may assume that $X$ is an affine scheme. Put $X = \Spec B$. If we show the following claim, then we get the assertion by taking $Z = \pi^{-1}(V(I))$ because a constructible subset which is closed under generalization is open.
    \begin{claim}\label{claim:claim1}
        Assume that $Z \subset V(IB)$ is a closed subset. Then $\mathsf{CM}(X) \cap Z$ is constructible in $Z$.
    \end{claim}
    \begin{pfclaim}
        We proceed by Noetherian induction. Let us suppose that any proper closed subset of $Z$ satisfies the conclusion of the claim. We may assume that $Z = V(\mathfrak{q})$ for some prime ideal $\mathfrak{q} \in V(IB)$. If $\mathfrak{q} \notin \mathsf{CM}(X)$, then $\mathsf{CM}(X) \cap Z = \emptyset$. This is clearly constructible. Assume that $\mathfrak{q} \in \mathsf{CM}(X)$. Since $A/I$ is $\mathsf{CM}$-Q0, $B/\mathfrak{q}$, a finitely generated $A/I$-algebra, is also $\mathsf{CM}$-Q0 by Remark~\ref{rem:CM-quasi-excellent} (as $\mathsf{CM}$-Q0 is equivalent to $\mathsf{CM}$-2). Hence, $\mathsf{CM}(B/\mathfrak{q})$ contains a non-empty open subset of $\Spec B/\mathfrak{q}$. Therefore, we can take $f \in B \setminus \mathfrak{q}$ such that $V(\mathfrak{q}) \cap D(f) \subset \mathsf{CM}(B) = \mathsf{CM}(X)$ by Lemma~\ref{lem:lem1-1}. Hence, we have
        \begin{align*}
            \mathsf{CM}(X) \cap (V(\mathfrak{q}) \setminus V(\mathfrak{q} + fB)) = V(\mathfrak{q}) \setminus V(\mathfrak{q} + fB).
        \end{align*}
        This set is open in $Z$. We have the following decomposition:
        \begin{align*}
            \mathsf{CM}(X) \cap Z = \left( \mathsf{CM}(X) \cap V(\mathfrak{q} + fB) \right) \cup \left( \mathsf{CM}(X) \cap (V(\mathfrak{q}) \setminus V(\mathfrak{q} + fB)) \right). 
        \end{align*}
        The former set is constructible in $Z$ by the induction hypothesis. The latter set is open in $Z$. Hence, $\mathsf{CM}(X) \cap Z$ is constructible in $Z$. This completes the proof of the claim.
    \end{pfclaim}
    This completes the proof of Lemma~\ref{lem:lem1-2}.
\end{proof}

The following lemma is the key step in which we use \v{C}esnavi\v{c}ius' Macaulayfication theorem. This is a Cohen--Macaulay analogue of \cite[Lemma~4.4]{KS21}.

\begin{lemma}\label{lem:lem1-3}
    Let $B$ be a Noetherian domain, $\mathfrak{q} \in \Spec B$ and $J \subset B$ an ideal. Assume that $B_\mathfrak{q}$ is a $\mathsf{CM}$\textnormal{-G-}ring and $B/J$ is $\mathsf{CM}$\textnormal{-Q0}. Then, there exists $b \in B \setminus \mathfrak{q}$, a scheme $X_b$ and a projective and surjective morphism $\phi_b \colon X_b \to \Spec B[b^{-1}]$ such that
    \begin{align*}
        \phi_b^{-1}\left( \Spec \frac{B[b^{-1}]}{JB[b^{-1}]} \right) \subset \mathsf{CM}(X_b).
    \end{align*}
\end{lemma}
\begin{proof}
    Since $B_\mathfrak{q}$ is $\mathsf{CM}$-quasi-excellent by Proposition \ref{prop:semilocalCMG}, there exists a Macaulayfication 
    \begin{align}\label{eq:macaulay}
        \pi_\mathfrak{q} \colon \tilde{X}_\mathfrak{q} \to \Spec B_\mathfrak{q}
    \end{align}
    by \cite[Theorem 1.6]{Ces21}. Since $\pi_\mathfrak{q}$ is projective, there exists $b_0 \in B \setminus \mathfrak{q}$ and a projective morphism
    \begin{align*}
        \pi_{b_0} \colon \tilde{X}_{b_0} \to \Spec B[b_0^{-1}]
    \end{align*}
    such that 
    \begin{align*}
        \tilde{X}_{b_0} \times_{\Spec B[b_0^{-1}]} \Spec B_\mathfrak{q} \cong \tilde{X}_\mathfrak{q}.
    \end{align*}
    Moreover, \eqref{eq:macaulay} is surjective and 
    \begin{align*}
        \Spec B_\mathfrak{q} = \varprojlim_{b \in B \setminus \mathfrak{q}} \Spec B[b^{-1}] = \varprojlim_{\substack{b \in B \setminus \mathfrak{q} \\ D(b) \subset D(b_0)}} \Spec B[b^{-1}].
    \end{align*}
    Hence, all the hypotheses of \cite[\href{https://stacks.math.columbia.edu/tag/081D}{Tag 081D}]{stacks-project} are satisfied. Applying \cite[\href{https://stacks.math.columbia.edu/tag/07RR}{Tag 07RR}] {stacks-project}, we obtain an element $\tilde{b}\in B\setminus\mathfrak{q}$ such that
    \begin{align*}
        \pi_{\tilde{b}} \colon X_{\tilde{b}} \to \Spec B[\tilde{b}^{-1}]
    \end{align*}
    is projective and surjective, where $X_{\tilde{b}} \coloneqq \tilde{X}_{b_0} \times_{\Spec B[b_0^{-1}]} \Spec B[\tilde{b}^{-1}]$. Set 
    \begin{align*}
        Z_{\tilde{b}} \coloneqq \Spec \frac{B[\tilde{b}^{-1}]}{JB[\tilde{b}^{-1}]}, \qquad W_{\tilde{b}} \coloneqq \pi_{\tilde{b}}^{-1}(Z_{\tilde{b}}) \setminus \mathsf{CM}(X_{\tilde{b}}). 
    \end{align*}
    Since $B/J$ is $\mathsf{CM}$-Q0, $B[\tilde{b}^{-1}]/JB[\tilde{b}^{-1}]$ is also $\mathsf{CM}$-Q0. By Lemma~\ref{lem:lem1-2}, we have that $\mathsf{CM}(X_{\tilde{b}}) \cap \pi_{\tilde{b}}^{-1}(Z_{\tilde{b}})$ is open in $\pi_{\tilde{b}}^{-1}(Z_{\tilde{b}})$. Hence, $W_{\tilde{b}}$ is closed in $\pi_{\tilde{b}}^{-1}(Z_{\tilde{b}})$. If $W_{\tilde{b}} = \emptyset$, then $b = \tilde{b}$ works. Assume that $W_{\tilde{b}} \neq \emptyset$. Since $\pi_{\tilde{b}}$ is projective, $\pi_{\tilde{b}}$ is proper. Hence, $\pi_{\tilde{b}}(W_{\tilde{b}})$ is closed in $\Spec B[\tilde{b}^{-1}]$. If $\mathfrak{q} \in \pi_{\tilde{b}}(W_{\tilde{b}})$, then we can find $s \in W_{\tilde{b}}$ such that $\pi_{\tilde{b}}(s) = \mathfrak{q}$. Moreover, we can find $\tilde{s} \in X_{\tilde{b}} \times_{\Spec B[\tilde{b}^{-1}]} \Spec B_\mathfrak{q} \cong \tilde{X}_\mathfrak{q}$ such that $\tilde{s} \mapsto s$ on the base change $X_{\tilde{b}} \times_{\Spec B[\tilde{b}^{-1}]} \Spec B_\mathfrak{q} \to X_{\tilde{b}}$. Hence, we have that $\mathcal{O}_{X_{\tilde{b}}, s} \cong \mathcal{O}_{\tilde{X}_\mathfrak{q}, \tilde{s}}$ is Cohen--Macaulay. This contradicts the fact that $s \in W_{\tilde{b}}$. Therefore, we have $\mathfrak{q} \notin \pi_{\tilde{b}}(W_{\tilde{b}})$. We can take $c \in B \setminus \mathfrak{q}$ such that $\Spec B[(\tilde{b}c)^{-1}] \cap \pi_{\tilde{b}}(W_{\tilde{b}}) = \emptyset$. Let $b = \tilde{b}c$, $X_b = X_{\tilde{b}} \times_{\Spec B[\tilde{b}^{-1}]} \Spec B[b^{-1}]$ and $\phi_b \colon X_b \to \Spec B[b^{-1}]$ be the induced morphism. 
    \[\begin{tikzcd}
        {X_b} & {\Spec B[b^{-1}]} \\
        {X_{\tilde{b}}} & {\Spec B[\tilde{b}^{-1}]}
        \arrow["{\phi_b}", from=1-1, to=1-2]
        \arrow[from=1-1, to=2-1]
        \arrow[from=1-2, to=2-2]
        \arrow["{\pi_{\tilde{b}}}"', from=2-1, to=2-2]
    \end{tikzcd}\]
    $\phi_b$ is also projective and surjective since $\pi_{\tilde{b}}$ is projective and surjective. Moreover, we have $W_{\tilde{b}} \cap X_b = \emptyset$. Therefore, we have
    \begin{align*}
        \phi_b^{-1}\left( \Spec \frac{B[b^{-1}]}{JB[b^{-1}]} \right) \subset \mathsf{CM}(X_b).
    \end{align*}
\end{proof}

\begin{remark}\label{rem:empty}
    The inverse image appearing in Lemma~\ref{lem:lem1-3} may be empty. Indeed, if $\sqrt{J}\not\subseteq\mathfrak q$, then one may take $b\in\sqrt{J}\setminus\mathfrak q$, in which case $B[b^{-1}]/JB[b^{-1}]$ is the zero ring, and the assertion becomes trivial. For example, this occurs when $J=B$. On the other hand, if $\sqrt{J}\subseteq\mathfrak q$, then every $b\in B\setminus\mathfrak q$ satisfies $b\notin\sqrt{J}$. Hence, $\operatorname{Spec}(B[b^{-1}]/JB[b^{-1}])$ is nonempty, and so is its inverse image under the surjective morphism $\phi_b$. Thus, it is the latter case that constitutes the nontrivial content of the lemma.
\end{remark}

Now, we are ready to prove Theorem~\ref{thm:main1}.

\begin{proof}[Proof of Theorem~\ref{thm:main1}]
    Taking $\mathfrak{p} \in \Spec A$ arbitrarily, we show that $\mathsf{CM}(B)$ contains a nonempty open subset of $\Spec B$, where $B \coloneqq A/\mathfrak{p}$. Since $A$ is a $\mathsf{CM}$\textnormal{-G-}ring, $B_\mathfrak{q}$ is $\mathsf{CM}$-quasi-excellent for any $\mathfrak{q} \in \Spec B$. By Lemma~\ref{lem:lem1-3}, we can take $b_{\mathfrak{q}} \in B \setminus \mathfrak{q}$, a scheme $X_{b_{\mathfrak{q}}}$ and a projective and surjective morphism $\phi_{b_{\mathfrak{q}}} \colon X_{b_{\mathfrak{q}}} \to \Spec B[b_{\mathfrak{q}}^{-1}]$ such that
    \begin{align*}
        \phi_{b_{\mathfrak{q}}}^{-1}\left( \Spec \frac{B[b_{\mathfrak{q}}^{-1}]}{IB[b_{\mathfrak{q}}^{-1}]} \right) \subset \mathsf{CM}(X_{b_{\mathfrak{q}}}). 
    \end{align*}
    Hence, we have
    \begin{align*}
        \Spec B = \bigcup_{\mathfrak{q} \in \Spec B} \Spec B[b_{\mathfrak{q}}^{-1}].
    \end{align*}
    Because $\Spec B$ is quasi-compact, we can take $b_1, \ldots, b_n \in \left\{ b_{\mathfrak{q}} \mid \mathfrak{q} \in \Spec B \right\}$ such that
    \begin{align}\label{eq:decomp}
        \Spec B = \bigcup_{i=1}^n \Spec B[b_i^{-1}].
    \end{align}
    Note that $\phi_{b_i} \colon X_{b_i} \to \Spec B[b_i^{-1}]$ is projective and surjective for each $i = 1, \ldots, n$ and 
    \begin{align*}
        \phi_{b_i}^{-1}\left( \Spec \frac{B[b_i^{-1}]}{IB[b_i^{-1}]} \right) \subset \mathsf{CM}(X_{b_i}). 
    \end{align*}
    By generic flatness, we can take $0 \neq c \in B$ such that 
    \begin{align}\label{eq:surj}
        \phi_{b_i}^{-1}\left( \Spec B[(cb_i)^{-1}] \right) \to \Spec B[(cb_i)^{-1}]
    \end{align}
    is flat for each $i = 1, \ldots, n$. We show that $B[c^{-1}]$ is Cohen--Macaulay. Let $P \in \Spec B[c^{-1}]$. We have to show that $B_P$ is Cohen--Macaulay. We take a maximal ideal $\mathfrak{m}$ of $B$ with $P \subset \mathfrak{m}$. By the decomposition \eqref{eq:decomp}, we can take $i \in \{ 1, \ldots, n \}$ such that $\mathfrak{m} \in \Spec B[b_i^{-1}]$. Hence, we have $P, \mathfrak{m} \in \Spec B[b_i^{-1}]$. Since \eqref{eq:surj} is surjective by Lemma \ref{lem:lem1-3}, we can take $x_P \in X_{b_i}$ such that $\phi_{b_i}(x_P) = P$. Since $\phi_{b_i}$ is continuous, we have $\phi_{b_i}(\overline{\left\{ x_P \right\}}) \subset \overline{\left\{ P \right\}}$. Clearly, $P \in \phi_{b_i}(\overline{\left\{ x_P \right\}})$. We have
    \begin{align*}
        \overline{\left\{ P \right\}} \subset \overline{\phi_{b_i}(\overline{\left\{x_P \right\}})} \subset \overline{\left\{ P \right\}} = V(P).
    \end{align*}
    Since $\phi_{b_i}$ is proper, $\phi_{b_i}(\overline{\left\{ x_P \right\}})$ is closed. $P \subset \mathfrak{m}$ means $\mathfrak{m} \in \overline{\left\{ P \right\}} = \phi_{b_i}(\overline{\left\{ x_P \right\}})$. Hence, there exists $x_\mathfrak{m} \in \overline{\left\{ x_P \right\}}$ such that $\phi_{b_i}(x_\mathfrak{m}) = \mathfrak{m}$. Since $I \subset \rad(A)$, $IB \subset \mathfrak{m}$. This means that $\mathfrak{m} \in \Spec B[b_i^{-1}]/IB[b_i^{-1}]$. We find that 
    \begin{align*}
        x_\mathfrak{m} \in \phi_{b_i}^{-1}\left( \Spec \frac{B[b_i^{-1}]}{IB[b_i^{-1}]} \right) \subset \mathsf{CM}(X_{b_i}).
    \end{align*}
    Therefore, $\mathcal{O}_{X_{b_i}, x_\mathfrak{m}}$ is Cohen--Macaulay. Hence, so is $\mathcal{O}_{X_{b_i}, x_P}$. The induced morphism
    \begin{align*}
        B_P = \mathcal{O}_{\Spec B[(cb_i)^{-1}], P} \to \mathcal{O}_{X_{b_i}, x_P}
    \end{align*}
    is faithfully flat. Therefore, $B_P$ is Cohen--Macaulay. This completes the proof of Theorem~\ref{thm:main1}.
\end{proof}

\bibliographystyle{amsalpha-inits}
\bibliography{references}

\end{document}